\documentclass[11pt]{article}

\usepackage{amsfonts}
\begin{document}

\baselineskip = 18 pt

\parskip = \the\baselineskip


\parindent = 0 pt

\newtheorem{thm}{Theorem}[section]

\newtheorem{lem}[thm]{Lemma} \newtheorem{definition}[thm]{Definition}

\newtheorem{pre-note}[thm]{Note}
\newenvironment{note}{\begin{pre-note}\rm}{\end{pre-note}}

\newtheorem{pre-proof}[thm]{Proof}
\newenvironment{proof}{\begin{pre-proof}\rm}{$\quad\bullet$\end{pre-proof}}

\newtheorem{pre-example}[thm]{Example}
\newenvironment{example}{\begin{pre-example}\rm}{\end{pre-example}}

\newtheorem{pre-numeq}[thm]{(}
\newenvironment{numeq}{\begin{pre-numeq}\rm\bf)$\quad\displaystyle}
{$\end{pre-numeq}}


\newcommand{\Aa}{\mathbb{A}} \newcommand{\RR}{\mathbb{R}}
\newcommand{\FF}{\mathbb{F}} \newcommand{\CC}{\mathbb{C}}
\newcommand{\QQ}{\mathbb{Q}} \newcommand{\HH}{\mathbb{H}}
\newcommand{\ZZ}{\mathbb{Z}} \newcommand{\NN}{\mathbb{N}}

\newcommand{\Cc}{{\cal C}} \newcommand{\Ss}{{\cal S}}

\newcommand{\Real}{\mathop{\rm Re}}
\newcommand{\inn}{\mathop{_{\rm \! inn}}}
\newcommand{\out}{\mathop{_{\rm \! out}}}

\newcommand{\LLim}{\mathop{\rm l\cdot i\cdot m}}

\begin{flushright} {\small 
math/0001013
to appear in the Journal of Number Theory} \end{flushright}

\vskip 1cm

\begin{center}
{\Large An adelic causality problem related to abelian
L--functions}\\
\vskip 1 cm
Jean-Fran\c{c}ois Burnol\\
January 2000, revised July 2000\\

\end{center}

\vfill {\bf Abstract:} We associate to the global field $K$ a
Lax--Phillips scattering which has the property of causality if and
only if the Riemann Hypothesis holds for all the abelian L--functions
of $K$. As a Hilbert space closure problem this provides an adelic
variation on a theme initiated by Nyman and Beurling. The adelic
aspects are related to previous work by Tate, Iwasawa and Connes.  
\vfill

Universit\'e de Nice\--\ Sophia Antipolis\\
Laboratoire J.-A. Dieudonn\'e\\
Parc Valrose\\
F-06108 Nice C\'edex 02\\
France\\
burnol@math.unice.fr\\

\clearpage 

\section{Adeles, Ideles and Zeros}

We express the Riemann Hypothesis for abelian L--functions as a
Hilbert space closure property (theorem {\bf\ref{AdelicNB}} below).
This takes place within the adelic set-up used by Tate \cite{Tat50}
(1950) and Iwasawa \cite{Iwa52} (1952) to establish the functional
equations of these L--functions. We treat simultaneously the number
field and function field cases (the Tate--Iwasawa ideas have been
adapted by Weil to the function field case in \cite{Wei67}). Our
approach is Hilbert space-theoretical. We take our hint from Nyman's
equivalent formulation of the original Riemann Hypothesis \cite{Nym50}
(1950). Beurling \cite{Beu49} (1949, for the disc), and  Lax
\cite{Lax59} (1959, for the half-plane), described the invariant
subspaces of the Hardy spaces and, as is explained in \cite{Bal98}
(see also \cite{Ber84} and \cite{Ber98} for Beurling's
$L^p$--extension \cite{Beu55} (1955)), this description is the
conceptual element behind Nyman's thorem. We devote a section to
explain (without mention of adeles and ideles) what our construction
amounts to for the Riemann zeta function. It is technically of a very
straightforward nature, its only deeper aspects being embedded in the
Beurling--Lax theory.

We associate to the global field an adelic Lax--Phillips scattering
\cite{Lax67} (1967). All axioms (where the idele class group replaces
the more usual $\ZZ$ or $\RR$) are satisfied, except possibly the
\emph{causality axiom} which we show to be equivalent to the Riemann
Hypothesis (this is our main result, theorem {\bf\ref{Causality}}).
The validity of one of the axioms is related to an observation of
Connes \cite[proof of VIII.5]{Con99}. The study of connections between
the Riemann zeta function and scattering theory is at least thirty
years old. In particular the Faddeev--Pavlov study of scattering for
automorphic functions \cite{Fad72} (1972), further developped by Lax
and Phillips in their book \cite{Lax76} (1976), has attracted
widespread attention. In their approach the scattering matrix is
directly related to the values taken by the Riemann zeta function on
the line $\Real(s) = 1$, and the Riemann Hypothesis itself is
equivalent to some decay properties of scattering waves. Another
well-known instance is the approach of De Branges (\cite{Bra86} 1986,
\cite{Bra94} 1994) within the theory of Hilbert spaces of entire
functions, also related to scattering (Conrey and Li have recently
pointed out some difficulties of this approach (\cite{Conr98}, 1998)).
The connection between our scattering process and the Riemann zeta
function (or more generally an abelian L--function) is the following:
each `bad' zero ($\Real(\rho) > {1\over2}$) appears as a pole of the
scattering operator, where there should be none, if the process was
causal. But if the Riemann Hypothesis holds, then the scattering
itself is of a trivial nature, and says absolutely nothing on the
zeros on the critical line. We point out that the same holds with the
positivity criterion of Weil (\cite{Wei52} 1952, \cite{Wei72} 1972):
the Weil distribution if of positive type if and only if the Riemann
Hypothesis holds, but beyond that, positivity tells nothing on the
location of the zeros except that they are indeed on the critical
line. Our formulation applies equally well to function fields and
number fields: this is as in Weil's positivity approach (especially
when formulated as in \cite{Bur99b}), and as in the work of Connes
(\cite{Con99}, 1999). The infinite places cause us less trouble than
in \cite{Wei52} and \cite{Con99}. Our sole motivation in formulating
the Riemann Hypothesis in a novel manner is the hope that creators of
other tools, of a deeper nature than those used here, would
incorporate the gained insight in their design constraints. An obvious
deficiency of this paper is its inability to achieve an alternative
proof of the Riemann Hypothesis in the function field case, where it
is not an hypothesis but a well-known theorem.

Let $K$ be a global field (an {\bf A}--field in the terminology of
Weil \cite{Wei67}): either an algebraic number field or a field
finitely generated and of transcendence degree $1$ over a finite
field. We briefly review some normalizations. The adele ring $\Aa_K$
is its own Pontrjagin dual. The set of characters (additive, unitary)
for which $K$ (diagonally embedded) is its own annihilator is
non-empty (and a single orbit under the action of $K^\times$). We pick
one such good character and let the additive Fourier transform ${\cal
F}$ be defined with respect to it (and the corresponding self--dual
Haar measure, which is in fact independent of the choice made). On
each local multiplicative group $K_\nu^\times$ we write $d^*v_\nu$ for
the multiplicative measure which assigns volume $1$ to the units
(finite place) or is ${dx\over2|x|}$ (real place) or
$drd\theta\over\pi r$ (complex place). On the idele group
$\Aa_K^\times$ (also seen as a subset of $\Aa_K$) we use $d^*v =
\prod_\nu d^*v_\nu$, and on the idele class group $\Cc_K =
\Aa_K^\times/K^\times$ we use the Haar measure $d^*u$ which (function
field case) assigns volume $1$ to the units or (number field case) is
pushed down to ${dt\over t}$ under $t = |u| = |v| = \prod_\nu
|v_\nu|_\nu$ ($v\in \Aa_K^\times, u = \overline{v}$).

Let $\Ss(\Aa_K)$ be the vector space of Bruhat--Schwartz functions.

\begin{definition} \begin{eqnarray*} E:  \Ss(\Aa_K) &\rightarrow&
(\Cc_K \rightarrow \CC) \\ \varphi(x)   &\mapsto&
f(\overline{v})=\sqrt{|v|}\sum_{q\in K^\times} \varphi(qv)  -
{\int_{\Aa_K} \varphi(x)\, dx \over \sqrt{|v|}}\\ \end{eqnarray*}
\end{definition}

For functions satisfying the additional conditions $\varphi(0) =
\int_{\Aa_K} \varphi(x)\, dx = 0$, $E$ is a tool at the heart of the
constructions of Connes in \cite{Con99}. For technical, class-field
theoretical, reasons, we do not impose any vanishing condition. The
map $E$ is related to the ideas of Tate \cite{Tat50} and Iwasawa
\cite{Iwa52}, and is especially tuned for Hilbert space matters, as
expressed in the following lemma:

\begin{lem}\label{ELemma} $E(\Ss(\Aa_K)) \subset L^2(\Cc_K, d^*u)$ and
is dense in it. The Fourier--Mellin transform of $E(\varphi)$, as a
function of the unitary characters of $\Cc_K$, is, up to a
multiplicative constant depending only on $K$, equal to the Tate
L--functions associated to $\varphi$ (restricted to the critical
line). \end{lem}

\begin{note} As has already been noted by Connes \cite[proof of
VIII.5]{Con99}, $E(\Ss_{00})$ is dense in $L^2(\Cc_K, d^*u)$, where
$\Ss_{00} = \{ \varphi \in \Ss(\Aa_K)\ |\ \varphi(0) = \int_{\Aa_K}
\varphi(x)\, dx = 0 \}$. \end{note}

The idele group acts on $\Ss(\Aa_K)$ ($U(v)\cdot \varphi(x) = {1 \over
\sqrt{|v|}}\varphi({x\over v})$) and on $L^2(\Cc_K, d^*u)$ ($U(v)\cdot
f(u) = f({u\over {\overline{v}}})$), and $E$ intertwines the two
actions. Furthermore the Poisson-Tate summation formula shows that $E$
intertwines the Fourier transform ${\cal F}$ on $\Aa_K$ with the
inversion $I$ ($f(u) \mapsto f({1\over u})$) on $\Cc_K$. Each idele
$v$ defines an adelic parallelepiped $$P(v) = \left\{ x =
(x_\nu)\in\Aa_K\ |\ \forall\ \nu\ |x_\nu|_\nu \leq
|v_\nu|_\nu\right\}$$ whose volume is proportional to $|v|$.

\begin{definition} $$\Ss_{\leq1} = \{ \varphi \in \Ss(\Aa_K)\ |\
\exists v\in \Aa_K^\times: |v| = 1\hbox{\rm\ and supp}(\varphi)\subset
P(v)\}$$ $$\widetilde{\Ss_{\leq1}} = \{ \varphi \in \Ss(\Aa_K)\ |\
\exists v\in \Aa_K^\times: |v| = 1\hbox{\rm\ and supp}({\cal
F}(\varphi))\subset P(v)\}$$ $${\cal D}_+ = {E(\Ss_{\leq 1})}^\perp$$
$${\cal D}_- = {E(\widetilde{\Ss_{\leq1}})}^\perp$$ \end{definition}

\begin{lem} The Lax--Phillips scattering axioms (\cite{Lax67}, with
$\ZZ$ or $\RR$ replaced with $\Cc_K$) are satisfied for the
``incoming'' subspace ${\cal D}_-$ $$|\lambda|\leq 1 \Rightarrow
U(\lambda){\cal D}_-\subset{\cal D}_-\quad\bigwedge U(\lambda){\cal
D}_- = \{0\}\quad\overline{\bigvee U(\lambda){\cal D}_-} =
L^2(\Cc_K,d^*u)$$ and for the ``outgoing'' subspace ${\cal D}_+$
$$|\lambda|\geq 1 \Rightarrow U(\lambda){\cal D}_+\subset{\cal
D}_+\quad\bigwedge U(\lambda){\cal D}_+ = \{0\}\quad\overline{\bigvee
U(\lambda){\cal D}_+} = L^2(\Cc_K,d^*u)$$ \end{lem}

\begin{note} The property $\bigwedge U(\lambda){\cal D}_+ = \{0\}$ is
cousin to the density property $\overline{E(\Ss_{00})} = L^2(\Cc_K,
d^*u)$ noted by Connes. The property $\overline{\bigvee
U(\lambda){\cal D}_+} = L^2(\Cc_K,d^*u)$ is an easy corollary of the
Artin--Whaples product formula. As ${\cal D}_- = I({\cal D}_+)$ and as
$I$ is an isometry which interchanges dilations and contractions, the
axioms for ${\cal D}_-$ and ${\cal D}_+$ are equivalent. \end{note}

Our main result is:

\begin{thm}[A causality criterion]\label{Causality} The Riemann
Hypothesis holds for all abelian L--functions of $K$ if and only if
${\cal D}_-\perp{\cal D}_+$. \end{thm}

We also express the Riemann Hypothesis as a closure property. We need
a slightly technical definition first:

\begin{definition}\label{Aoperator} Let $A$ be the convolution
operator $$(A\cdot f)(u_0) = \int_{\Cc_K} a({u_0\over u})f(u)\, d^*u$$
where, in the number field case $$a(w) = \sqrt{|w|}\cdot{\bf
1}_{|w|\leq1}$$ and in the function field case ($q$ the cardinality of
the field of constants) $$a(w) =
\sqrt{|w|}\cdot(\sqrt{q}-{1\over\sqrt{q}})\cdot{\bf 1}_{|w|<1} + (1 -
{1\over\sqrt{q}}){\bf 1}_{|w|=1}$$ \end{definition}

\begin{definition} $$\HH^2 = \left\{ f \in L^2(\Cc_K,d^*u)\;|\
\mbox{\rm ess-supp}(f)\subset\{|u|\leq 1\} \right\}$$ \end{definition}

\begin{lem} The operator $V = 1 - A$ is a unitary operator on
$L^2(\Cc_K, d^*u)$, commuting with the regular action of $\Cc_K$, and
sending $\HH^2$ to (a subspace of) itself. \end{lem}

\begin{thm}[A closure criterion]\label{AdelicNB}
$V(\overline{E(\Ss_{\leq 1})}) \subset \HH^2$ with equality if and
only if the Riemann Hypothesis holds for all abelian L--functions of
$K$. \end{thm}

\section{The criterion for the Riemann zeta function}

When considering only the Riemann zeta function, Theorem
{\bf\ref{AdelicNB}} boils down to a variant of Nyman's criterion
\cite{Nym50}. Let us recall this criterion (see also \cite{Bal98},
\cite{Ber84}, \cite{Ber98}):

Let $\rho_\alpha(u) = \left\{\alpha\over u\right\} -
\alpha\left\{1\over u\right\}$, for $0<\alpha<1$, and $u\in (0,1)$
(with $\{\cdot\}$ the fractional part). Let $N$ be the closed span in
$L^2((0,1), du)$ of the functions $\rho_\alpha$. We also consider both
$N$ and $L^2((0,1), du)$ as closed subspaces of $L^2((0,\infty), du)$.

\begin{thm}[Nyman, 1950 \cite{Nym50}]\label{Nyman} The constant
function {\bf 1} on $(0,1)$ belongs to $N$ if and only if the Riemann
Hypothesis holds. \end{thm}

Note that $N$ is invariant under the semi-group of unitary
contractions $U(\lambda):f(u)\mapsto
\sqrt{{1\over\lambda}}f({u\over\lambda})$, $\lambda \leq 1$, $u > 0$
as $U(\lambda)\cdot\rho_\alpha =
\sqrt{{1\over\lambda}}(\rho_{\alpha\lambda} - \alpha\rho_\lambda)$.
So, it will contain the constant function ${\bf 1}$ (hence all step
functions) if and only if it actually coincides with all of
$L^2((0,1), du)$.

For the proof one considers the Mellin transform: $$f(u) \mapsto
\widehat{f}(s) = \int_0^1 f(u) u^{s-1} du $$ which by a Paley--Wiener
theorem establishes an isometry between $L^2((0,1),du)$ and the Hardy
space $\HH^2(\Real(s)>{1\over2})$ of analytic functions with bounded
norm $$\Vert F \Vert^2 = \sup_{\sigma>{1\over2}}\int_{\Real(s) =
\sigma} |F(s)|^2 \,{|ds|\over 2\pi}$$ Such functions $\widehat{f}(s)$
have (a.e.) boundary values also obtained as $$\widehat{f}({1\over2} +
i\tau) = \LLim_{\epsilon\to0}\int_\epsilon^1 \sqrt{u}f(u)
u^{i\tau}{du\over u}$$ equivalently as the Fourier--Plancherel
transform of $e^{t/2}f(e^t)$, $t\leq0$.

The unitary semi-group considered above acts on
$\HH^2(\Real(s)>{1\over2})$ as $F(s) \mapsto \lambda^{s-{1\over2}}
F(s)$, and Lax \cite{Lax59} has described the closed subspaces
invariant under this action. It can be directly shown (see
\cite{Hof62}) that the conformal representation $$w = {s-1 \over s}$$
$$g(w) = s\cdot F(s)$$ establishes an isometry between
$\HH^2(\Real(s)>{1\over2})$ and $\HH^2(|w|< 1)$ which identifies the
invariant subspaces of the former with closed subspaces invariant
under the shift $g(w) \mapsto w\cdot g(w)$ for the latter. These were
described by Beurling \cite{Beu49} and we learn that the
``continuous'' case (Lax) and ``discrete'' case (Beurling) are
completely equivalent (this equivalence is also a corollary to the
conformal invariance of Brownian motion on the complex numbers).

The Beurling--Lax recipe to determine an invariant closed subspace
such as $N$ is to look at the Mellin transforms of the functions
$\rho_\alpha(u)$'s: $$\widehat{\rho_\alpha}(s) =
{{\alpha-\alpha^s}\over s}\zeta(s)$$ and at the ``greatest lower bound
of their inner factors'': first there will be the Blaschke product
$$B(s) = \prod_{\zeta(\rho) = 0,\ \Real(\rho)>{1\over2}}{s - \rho
\over s - (1 -\bar{\rho})}{1 - \bar{\rho}\over\rho}\left|{\rho\over
1-\rho}\right|$$ where the zeros appear according to their
multiplicities, then an inner factor associated to a singular measure
on the critical line (the analytic continuation of $\zeta(s)$ implies
its non-existence), and a final inner factor $\lambda^{s-{1\over2}}$
($0<\lambda\leq1$). We argue that $\lambda = 1$ as follows:
$\lambda^{s-{1\over2}}\HH^2$ is the Mellin transform of
$L^2((0,\lambda),du)$ which contains $N$ only if $\lambda = 1$
(obviously). 

Bercovici and Foias \cite[2.1]{Ber84} prove $\lambda = 1$ in the
following manner: if
$\widehat{\rho_\alpha}(s)=\lambda^{s-{1\over2}}f(s)$ for some
$f(s)\in\HH^2(\Real(s)>{1\over2})$ then $\widehat{\rho_\alpha}(\sigma)
= O(\lambda^\sigma)$ for $\sigma\to+\infty$. Indeed
\footnote{I thank the referee for correcting my incomplete
understanding of the Bercovici--Foias proof at this point.}
$f(s)$ is $O(1)$ in any half-plane $\Real(s)\geq{1\over2} +
\varepsilon$, $\varepsilon >0$ (this follows from its Cauchy integral
representation or from $\widehat{f}(s) = \int_0^1 f(u) u^{s-1} du$ and
Cauchy-Schwarz). But obviously
$\lim_{\sigma\to\,+\infty}\,\sigma\cdot\widehat{\rho_\alpha}(\sigma)
\neq 0$, thus giving a contradiction if $\lambda < 1$ . The following
lemma, of independent interest, could also have been used:

\begin{lem}\label{GLemma} If $F(s) \in \HH^2$ is $O(|s|^K)$ on the
critical line, then its outer factor $F\out(s)$ is $O(|s|^K)$ on the
entire closed half-plane. \end{lem} \begin{proof} One has
$$\log(|s\,F\out(s)|) = \int_{Re(s_0) = {1\over2}}
\log(|s_0\,F(s_0)|){2\Real(s) -1 \over |s-s_0|^2}\, {|ds_0|\over
2\pi}$$ and $$\log(|s|) = \int_{Re(s_0) = {1\over2}}
\log(|s_0|){2\Real(s) -1 \over |s-s_0|^2}\, {|ds_0|\over 2\pi}$$ hence
the result \end{proof}

Let us add a few more words to this discussion of Nyman's theorem. As
$$\int_0^1 \left\{1\over u\right\}u^{s-1}\, du = {1\over s-1} -
{\zeta(s)\over s}$$ (for $\Real(s) > 0$) and ${s-1\over s}\cdot{1\over
s-1} = {1\over s} = \int_0^1 u^{s-1}\, du$ we see that ${s-1\over
s}{\zeta(s)\over s}$ belongs to $\HH^2(\Real(s)>{1\over2})$. The
unitary operator $V$ on $L^2((0, \infty), du)$ given by the multiplier
${s-1\over s}$ in the spectral representation acts as $$f(u) \mapsto
f(u) - \int_u^\infty {1\over t}f(t)\,dt$$ As ${\zeta(s)\over s} = -
\int_0^\infty \{{1\over u}\} u^{s-1}\,du$ (for $0< \Real(s) <1$) we
obtain after a straightforward computation: $${s-1\over
s}{\zeta(s)\over s} = \int_0^1 A(u) u^{s-1}\,du$$ $$A(u) = [{1\over
u}]\log(u) + \log([{1\over u}]!) + [{1\over u}]$$ Stirling's formula
implies $A(u) = {1\over2}\log({1\over u}) + O(1)$ so this integral
representation is valid for $\Real(s) > 0$. As $A(u) = 1 + \log(u)$
for ${1\over 2} < u\leq 1$ there is no inner factor of the type
$\lambda^{s-{1\over2}}$ with $\lambda < 1$. There is no other singular
factor thanks to the analytic continuation, so ${s-1\over
s}{\zeta(s)\over s}$ is the product of an outer factor with the
Blaschke product $B(s)$. Hence

\begin{thm} The Riemann Hypothesis holds if and only if ${s-1\over
s}{\zeta(s)\over s}$ is an outer function, or equivalently if the
functions $U(\lambda)\cdot A(u)$ ($0<\lambda\leq1$) span
$L^2((0,1),du)$. \end{thm}

The generalized Jensen's formula (see \cite{Hof62}) then implies a
formula first derived by Balazard, Saias and Yor:

\begin{thm}[\cite{Bal2}] $${1\over2\pi}\int_{\Real(s) = {1\over2}}
{\log|\zeta(s)|\over|s|^2}\,|ds| = \sum_{\zeta(\rho) = 0,\
\Real(\rho)>{1\over2}}\log\left|{\rho\over 1-\rho}\right|$$ \end{thm}

The only difference with the proof of Balazard, Saias and Yor is that
we do not need the general theory of Hardy spaces beyond that of
$\HH^2$, which is of a more elementary nature. This concludes our
discussion of Nyman's theorem. We now turn to some variations on this
theme (other variations have been considered by Bercovici and Foias in
\cite{Ber84} and \cite{Ber98}).

Let $\phi(x)$ be a smooth function on the real line with compact
support in $[0,1]$, and $\int_0^1 \phi(x)\,dx = 0$. The Mellin
transform $$\widehat{\phi}(s) = \int_0^\infty \phi(u) u^{s-1} du$$ is
an entire function, vanishing at $1$. We consider: $$T(\phi)(u) =
\sum_{n\geq 1} \phi(nu)\quad(u>0)$$ which is a smooth function of $u$
on $(0,\infty)$ with support in $(0,1]$. Its behavior when
$u\rightarrow 0$ is governed by the Poisson summation formula:
$$T(\phi)(u) = {1\over |u|} \sum_{n\in\ZZ} \psi({n\over u})$$ where
$\psi$ is the Fourier transform $\int \phi(y)e^{2\pi i\,xy}\,dy$ of
$\phi$ (hence belongs to the Schwartz space of rapidly decreasing
functions). So $$\forall K\quad T(\phi)(u) =_{u\rightarrow 0} O(u^K)$$
and the Mellin transform $$\widehat{T(\phi)}(s) = \int_0^1 T(\phi)(u)
u^{s-1} du$$ is an entire function. For $Re(s) > 1$
$$\widehat{T(\phi)}(s) = \zeta(s) \widehat{\phi}(s)$$ hence by
analytic continuation this holds true for all $s$.

Let $\Ss^0_{\leq1}$ be the vector space consisting of these functions
$\phi$, $\overline{\Ss^0_{\leq1}}$ its closure in $L^2((0,1), du)$ and
$K$ the closure of the vector space of functions $T(\phi)$. Both
$\overline{\Ss^0_{\leq1}}$ and $K$ are invariant under contractions,
hence described by the Beurling--Lax theory. One just has to take the
``greatest lower bound'' of the inner factors of the
$\widehat{\phi}(s)$'s (resp. the $\widehat{T(\phi)}(s)$'s). Obviously
$\overline{\Ss^0_{\leq1}}$ is the subspace perpendicular to the
constant {\bf 1} and this shows that the ``greatest lower bound'' for
the zeros of the $\widehat{\phi}(s)$'s is simply $s = 1$ with
multiplicity $1$. This cancels exactly the pole of the zeta function.
For the $\widehat{T(\phi)}(s)$'s the analytic continuation across the
critical line implies that the only possible singular factor is of the
type $\lambda^{s-{1\over2}}$ with $\lambda \leq 1$. For a suitably
chosen $\phi$, $T(\phi)$ does not vanish in $({1\over2},1)$ so
necessarily $\lambda = 1$. The conclusion is that $K$ coincides with
the space $N$ considered by Nyman. Thus:

\begin{thm} The Riemann Hypothesis holds if and only if the constant
function {\bf 1} belongs to the closure of $\{T(\varphi): \varphi\in
\Ss^0_{\leq1}\}\quad\bullet$ \end{thm}

We describe one more variation. Let $\Ss^{ev}$ be the vector space of
even Schwartz functions on $\RR$. Let, for $u > 0$: $$E(\varphi)(u) =
\sum_{n\geq 1} \varphi(nu) - {\int_0^\infty \varphi(x)\, dx \over u}$$
The Poisson summation formula gives $$E(\varphi)(u) = {1\over
u}\sum_{n\geq 1} {\cal F}(\varphi)({n\over u}) - {1\over2}\varphi(0)$$
so that $E(\varphi)(u)$ is $0(1)$ when $u \rightarrow 0$ and is
$O({1\over u})$ when $u \rightarrow \infty$ and belongs to $L^2(\RR_+,
du)$. Its Mellin transform $$\widehat{E(\varphi)}(s) = \int_0^\infty
E(\varphi)(u) u^{s-1} du$$ is absolutely convergent and analytic for
$0 < \Real(s) < 1$. It can be rewritten as $$\int_0^1 E(\varphi)(u)
u^{s-1} du  + \int_1^\infty \sum_{n\geq 1} \varphi(nu) u^{s-1} du +
{\int_0^\infty \varphi(x)\, dx\over s-1}$$ which is then valid in the
half-plane $\Real(s) > 0$. Then, for $\Real(s) > 1$, as $$
\int_0^\infty \sum_{n\geq 1} \varphi(nu) u^{s-1} du - \int_0^1
{\int_0^\infty \varphi(x)\, dx\over u} u^{s-1} du  + {\int_0^\infty
\varphi(x)\, dx\over s-1} $$ hence simply as $$\sum_{n\geq 1} n^{-s}\,
\int_0^\infty \varphi(u) u^{s-1} du = \zeta(s)\widehat{\varphi}(s)$$
which remains valid for $\Real(s) > 0$.

We now need to get rid of the pole of $\zeta(s)$ with the help of the
operator $V$ (which on $L^2(\RR_+,du)$ acts as ${s-1\over s}$ in the
spectral representation): $$V\cdot f(u) = f(u) - \int_u^\infty {1\over
v} f(v) dv$$ One checks $V\cdot {1\over u} = 0$ so $$VE(\varphi)(u) =
\sum_{n\geq 1} \varphi(nu) - \int_u^\infty \sum_{n\geq 1} \varphi(nv)
{dv\over v} = \sum_{n\geq 1} \varphi(nu) - \int_0^\infty [{v\over
u}]\varphi(v) {dv\over v}$$ Let $\Ss_{\leq1}$ be the vector space of
smooth even functions with support in $[-1,1]$. For
$\varphi\in\Ss_{\leq1}$, $VE(\varphi)$ has support in $(0,1]$ and its
Mellin transform ${s-1\over s}\zeta(s)\widehat{\varphi}(s)$ thus
belongs to $\HH^2$. As in the previous discussions, the Mellin
transform of the (closure of) $VE(\Ss_{\leq1})$ is the space of
multiples of the Blaschke product $B(s)\HH^2$. Hence:

\begin{thm} $\overline{VE(\Ss_{\leq1})} \subset \HH^2$ with equality
if and only if the Riemann Hypothesis holds. \end{thm}

Let $B$ be the unitary operator on $L^2(\RR_+,du)$ which acts in the
spectral representation as multiplication with $B(s)$. Let $${\cal
D}_+ = E(\Ss_{\leq1})^\perp = V^{-1}B\cdot(\HH^2)^\perp =
V^{-1}BI\cdot\HH^2$$ (where $I$ is the inversion $f(u) \mapsto {1\over
u}f({1\over u})$, or spectrally $s \mapsto 1-s$). Let $${\cal D}_- =
E({\cal F}(\Ss_{\leq1}))^\perp = I({\cal D}_+) = IV^{-1}BI\cdot\HH^2 =
VB^{-1}\cdot\HH^2$$ Then, in the terminology of Lax and Phillips
\cite{Lax67}, ${\cal D}_+$ (resp. ${\cal D}_-$) is an ``outgoing''
(resp. ``incoming'') space for the action of $\RR_+^\times$ on
$L^2(\RR_+,du)$. The scattering operator associated to them is $$S =
(V^{-1}B)^{-1}\cdot VB^{-1} = V^2B^{-2}$$ It is an invariant operator
whose spectral multiplier is $({s-1\over s})^2\cdot B(s)^{-2}$ and is
an inner function if and only if $B(s)$ has no zero in $\Real(s) >
{1\over2}$, that is if the Riemann Hypothesis holds. The scattering
multiplier is inner if and only if ${\cal D}_+ \perp {\cal D}_-$. So:

\begin{thm} $E(\Ss_{\leq1})^\perp \subset \overline{E({\cal
F}(\Ss_{\leq1}))}$ if and only if the Riemann Hypothesis holds.
\end{thm}

\section{An adelic scattering}

We now prove theorems {\bf \ref{Causality}} and {\bf \ref{AdelicNB}}.
Let $\Cc_K^1$ be the (compact) subgroup of idele classes of unit
modulus. There is some (non-canonical) isomorphism $\Cc_K = \Cc_K^1
\times N$, $N = \{ |u| : u\in \Cc_K \} \subset\RR_+^\times$. If $K$
has positive characteristic we let $q$ be the cardinality of the field
of constants. It is known that the module group $N$ is $q^\ZZ$. Each
character $\chi$ of $\Cc_K^1$ extends to a character of $\Cc_K$
trivial on $N$, which we still denote by $\chi$. At each place $\nu$
there is a local character $\chi_\nu$ from the embedding $K_\nu^\times
\rightarrow \Cc_K$. And $\chi$ is said to be ramified at $\nu$ if the
restriction of $\chi_\nu$ to the unit subgroup is non-trivial.

We start with the properties of \begin{eqnarray*} E:  \Ss(\Aa_K)
&\rightarrow& (\Cc_K \rightarrow \CC) \\ \varphi(x)   &\mapsto&
f(\overline{v})=\sqrt{|v|}\sum_{q\in K^\times} \varphi(qv)  -
{\int_{\Aa_K} \varphi(x)\, dx \over \sqrt{|v|}}\\ \end{eqnarray*} From
the definition one has $E(\varphi)(u) = O({1\over \sqrt{|u|}})$ when
$|u| \rightarrow \infty$, and as the Poisson-Tate formula gives
$$E\cdot{\cal F} = I\cdot E$$ one also has $E(\varphi)(u) =
O(\sqrt{|u|})$ when $|u| \rightarrow 0$. So indeed $$E(\Ss(\Aa_K))
\subset L^2(\Cc_K, d^*u)$$

Let $\chi$ be a unitary character on $\Cc_K$ (trivial on $N$). The
Fourier-Mellin transform (for $\Real(s) = {1\over2}$)
$$\widehat{E(\varphi)}(\chi,s) = \int_{\Cc_K}
E(\varphi(u))\chi(u)|u|^{s-{1\over2}}\,d^*u$$ is in fact absolutely
convergent and analytic for $0 < \Real(s) < 1$. It can be rewritten
(with $u = \overline{v}$, $v \in \Aa_K^\times$) as $$\int_{|u|\leq1}
E(\varphi)(u) \chi(u)|u|^{s-{1\over2}}\,d^*u  + \int_{|u|>1}
\sum_{q\in K^\times} \varphi(qv) \chi(u)|u|^{s}\,d^*u$$ $$-
\int_{\Aa_K} \varphi(x)\, dx \int_{|u|>1} \chi(u)|u|^{s-1}\,d^*u$$ The
integral $\int_{|u|>1} \chi(u)|u|^{s-1}\,d^*u$ (which vanishes if
$\chi \neq {\bf 1}$) is a meromorphic function $F_\chi(s)$, which can
be evaluated explicitely. One obtains (both in the number field and
function field cases) $$(\Real(s) > 1) \Rightarrow F_\chi(s) = -
\int_{|u|\leq1} \chi(u)|u|^{s-1}\,d^*u$$ So
$\widehat{E(\varphi)}(\chi,s)$ has a meromorphic continuation to
$\Real(s) > 0$ which, for $\Real(s) > 1$, coincides with
\medbreak
$$\int_{|u|\leq1} E(\varphi)(u) \chi(u)|u|^{s-{1\over2}}\,d^*u  +
\int_{|u|>1} \sum_{q\in K^\times} \varphi(qv) \chi(u)|u|^{s}\,d^*u$$
$$+ \int_{\Aa_K} \varphi(x)\, dx \int_{|u|\leq1}
\chi(u)|u|^{s-1}\,d^*u$$ $$= \int_{\Cc_K} \sum_{q\in K^\times}
\varphi(qv) \chi(u)|u|^{s}\,d^*u$$ $$= C(K) \int_{\Aa_K^\times}
\varphi(v) \chi(v)|v|^s d^*v$$ The constant $C(K)$ being as in Tate's
thesis \cite{Tat50} related to the way the measures $d^*u$ on $\Cc_K$
and $d^*v$ on $\Aa_K^\times$ differ. We recognize in the last integral
the Tate L--function $L(\varphi, \chi, s)$. The identity
$$\widehat{E(\varphi)}(\chi,s) = C(K) L(\varphi, \chi, s)$$ for
$\Real(s) = {1\over2}$ holds by analytic continuation. With this lemma
{\bf \ref{ELemma}} is proven.

We turn to the description of $\Delta = \overline{E(\Ss_{\leq1})}$.
The crucial thing is that it is invariant (obviously) under the
(unitary) action of the semi-group of contractions $\{|u| \leq 1\}$,
in particular under the action of the compact group $\Cc_K^1$. It thus
decomposes as a Hilbert space sum of isotypical components
$\Delta_\chi$, which we wish to compare to the isotypical components
of $\HH^2 = \left\{ f \in L^2(\Cc_K,d^*u)\;|\ \mbox{\rm
ess-supp}(f)\subset\{|u|\leq 1\} \right\}$. We do this in the spectral
representation using the Fourier--Mellin transform (in the function
field case we write $z = q^{-(s-{1\over2})}$).

Firstly it is a straightforward check that the $A$-operator
({\bf\ref{Aoperator}}) is an invariant operator whose action on $L^2$
is given by the following spectral multipliers $A(\chi,s)$: $$\chi
\neq 1 \Rightarrow A(\chi,s) = 0$$ $$A(1,s) = {1\over s}
\quad\mbox{\rm (number field case)}$$ $$A(1,z) = 1 - {1-\sqrt{q}z\over
\sqrt{q} -z} \quad\mbox{\rm (function field case)}$$ so that $V = 1 -
A$ is indeed a unitary (on $L^2$) invariant operator with multipliers
$$\chi \neq 1 \Rightarrow V(\chi,s) = 1$$ $$V(1,s) = {s-1\over s}
\quad\mbox{\rm (number field case)}$$ $$V(1,z) = {1-\sqrt{q}z\over
\sqrt{q} -z} \quad\mbox{\rm (function field case)}$$

From this spectral representation or with a direct computation we also
find the important identity $$V({1\over\sqrt{|u|}}\cdot{\bf 1}_{|u| >
1}) = -\alpha(K)\sqrt{|u|}\cdot{\bf 1}_{|u| \leq 1}$$ with $\alpha(K)
= 1$ (resp. $ {1\over\sqrt{q}}$) in the number field case (resp.
function field case). From the Artin--Whaples product formula we
obtain $E(\varphi)(u) = -{\int_{\Aa_K}\varphi(x)\,dx\over\sqrt{|u|}}$
for $|u| > 1$ and $\varphi \in \Ss_{\leq1}$. So we see that
$V(\Delta)$ is a subspace of $\HH^2$.  We now describe it exactly with
the help of the Beurling--Lax theory.

Let $S_f$ be the set of finite places of $K$, and $S_\infty$ the
(possibly empty) set of infinite places. Let $q_\nu$ be the
cardinality of the residue field at the finite place $\nu$, $\pi_\nu$
a uniformizer element of $K_\nu^\times$, which we also consider as an
element of $\Aa_K^\times$. The value $\chi(\pi_\nu)$ is independent of
the choice of $\pi_\nu$ if the character $\chi$ is unramified at
$\nu$. The (``incomplete'' in the number field case) L--function
associated to $\chi$ is $$L(\chi,s) = \prod_{\nu \in S_f,
\mbox{unramified}} {1\over {1 - \chi(\pi_\nu)q_\nu^{-s}}}$$

The Bruhat-Schwartz function $\varphi$ is built from local components,
all of them except finitely many being equal to the characteristic
function of the local integers, so its Tate L--function $L(\varphi,
\chi, s)$ is a multiple of $L(\chi,s)$ by a function holomorphic in
$\Real(s) > 0$. By lemma {\bf\ref{ELemma}} this implies that the
Paley--Wiener transform $\widehat{E(\varphi)}(\chi,s)$ ($\Real(s) >
{1\over2}$) vanishes at each bad zero with at least the same
multiplicity as $L(\chi,s)$.

\begin{definition} Let $B$ be the unitary invariant operator whose
spectral multiplier in the $\chi$-isotypical component of
$L^2(\Cc_K,d^*u)$ is the Blaschke product on the zeros (with
multiplicity) of the L--function $L(\chi,s)$ in the half-plane
$\Real(s) > {1\over2}$ (number field case) or the open disc $|z| < 1$
($z= q^{-(s-{1\over2})}$, function field case). \end{definition}

We will soon show that one can indeed build a convergent Blaschke
product with the bad zeros so that $B$ exists! (the function field
case is trivial as there are only finitely many). This being
temporarily admitted we have obtained $V(\Delta)\subset B\cdot\HH^2$.
And we prove

\begin{thm} $$V(\Delta) = B\cdot\HH^2$$ \end{thm}

We treat the function-field case first. We choose $\varphi_\nu$ to be
${\bf 1}_{|x|_\nu \leq 1}$ at a non-ramified place, and
$\overline{\chi_\nu(x)}\cdot{\bf 1}_{|x|_\nu = 1}$ at a ramified
place. With these choices we obtain $\varphi = \prod_\nu \varphi_\nu$
which belongs to $\Ss_{\leq 1}$ and for which (at first for $\Real(s)
> 1$): $$L(\varphi, \chi, s) = L(\chi,s)$$ We do not claim that
$E(\varphi)$ is $\chi$--equivariant, nevertheless this identity
combined with lemma {\bf\ref{ELemma}} and the inclusion $V(\Delta)
\subset \HH^2$ shows that $V(\chi,s)L(\chi,s)$ belongs to
$\HH^2(|z|<1)$. It is clear from the product representation that it
does not vanish at $z=0$, and it is known for $\chi = 1$ that the pole
at $s=1$ of the zeta function $Z_K(s)$ is of order $1$. Analytic
continuation across $|z| = 1$ implies the non-existence of a singular
inner factor. So the smallest closed subspace of $\HH^2(|z|<1)$
containing $V(\chi,s)L(\chi,s)$, and invariant under shifts, is
exactly $B(\chi,s)\HH^2$. The conclusion follows.

Let us now consider the case where $K$ is an algebraic number field. 
We define $\varphi_\nu(x_\nu)$ exactly as in the function field case
when $\nu$ is finite and as $\overline{\chi_\nu(x)}\cdot
g_\nu(|x|_\nu)$ at each infinite place, with $g_\nu$ a smooth function
on $\RR^\times_+$ with compact support in $(0,1)$. The product
function $\varphi(x) = \prod_\nu \varphi_\nu(x_\nu)$ then belongs to
$\Ss_{\leq 1}$ and $E(\varphi)$ has a Paley--Wiener transform
\begin{eqnarray*} & &\int_{\Cc_K}
E(\varphi)(u)\cdot\chi(u)|u|^{s-{1\over2}}\,d^*u \\ &=& C(K)
\int_{\Aa_K^\times} \varphi(v) \chi(v) |v|^s \, d^*v \\ &=& C(K)
L(\chi, s)\cdot\prod_{\nu\in S_\infty}\widehat{g_\nu}(s)\\
\end{eqnarray*} From this and the inclusion $V(\Delta) \subset \HH^2$
follows the existence of the Blaschke product $B(\chi,s)$ as promised
above. Furthermore it is clearly possible to choose the $g_\nu$ in
such a manner that $\widehat{g_\nu}(s)$ does not vanish at any $s$
prescribed in advance, and the existence of analytic continuation
accross the critical line then reduces the possibility of an inner
factor to $\lambda^{s-{1\over2}}$ with $\lambda \leq 1$. The
Bercovici--Foias argument implies as
in our discussion of Nyman's theorem that $\lambda = 1$. Finally it is
known that the pole of the zeta function ($\chi = 1$) has exact order
$1$. With all this the identity $V(\Delta) = B\cdot\HH^2$ is proven.
This completes the proof of the closure criterion {\bf\ref{AdelicNB}}.

Let ${\cal D}_+ = {E(\Ss_{\leq 1})}^\perp = \Delta^\perp =
V^{-1}B\cdot(\HH^2)^\perp$. Let $Z$ be the unitary operator which is
just $1$ in the number field case and $z$ (in each isotypical
component) in the function field case. Then $(\HH^2)^\perp =
Z^{-1}I\cdot\HH^2$ and ${\cal D}_+ = V^{-1}BZ^{-1}I\cdot\HH^2$. From
this follows $$\bigwedge U(\lambda){\cal D}_+ =
\{0\}\quad\overline{\bigvee U(\lambda){\cal D}_+} = L^2(\Cc_K,d^*u)$$
so that ${\cal D}_+$ indeed qualifies as an outgoing subspace and
${\cal D}_-$ as an incoming subspace. One has ${\cal D}_- =
IV^{-1}BZ^{-1}I\cdot\HH^2 = VB^{-1}Z\cdot\HH^2$. The Lax--Phillips
scattering operator associated to the pair $({\cal D}_+,{\cal D}_-)$
is an invariant unitary operator, unique up to a multiplicative
constant in each isotypical component. It is: $$S =
(V^{-1}B)^{-1}\cdot VB^{-1}Z = ZV^2B^{-2}$$ With the help of $S$ the
pair $({\cal D}_+,{\cal D}_-)$ is unitarily equivalent to
$((\HH^2)^\perp,S\cdot\HH^2)$. So it is an orthogonal pair if and only
if $S\cdot\HH^2\subset\HH^2$, if and only if $B = 1$, if and only if
the Riemann Hypothesis holds for all abelian L--functions of $K$. With
this the proof of the causality criterion {\bf \ref{Causality}} is
complete.

\begin{note} The reader of the monograph of Lax and Phillips
\cite[chapter 2]{Lax67} will perhaps be perplexed by the fact that
``causal'' means there ``inner with respect to the exterior domain
$|z| > 1$'' (in the discrete case). But this is because they represent
the semi-group leaving invariant the outgoing space with the help of
the non-negative powers of $z$. In our case we represent it with the
help of the non-negative powers of ${1\over z}$. So ``causal'' is to
be understood to mean ``inner with respect to the domain $|{1\over z}|
> 1$'' (that is $|z| < 1$). \end{note}

\begin{note} We have used $IBI = B^{-1}$. This follows from
$\overline{L(\chi,\overline{s})} = L(\overline{\chi},s)$ which implies
$B(\overline{\chi},\overline{s}) = \overline{B(\chi, s)}$
($=B(\chi,s)^{-1}$ for $\Real(s) = {1\over2}$). \end{note}

\begin{flushleft}
{\small Universit\'e de Nice\ --\ Sophia Antipolis\\
Laboratoire J.-A. Dieudonn\'e\\
Parc Valrose\\
F-06108 Nice C\'edex 02\\
France\\
burnol@math.unice.fr\\}
\end{flushleft}

\clearpage

\end{document}